\documentstyle[12pt,amstex,amssymb]{article}

\def\fin{e.\ g.\ } 
\def\this{i.e.\ }
\def\h{\hbar}
\def\p{\partial}
\def\w{\wedge}
\def\v{^{\vee }}

\def\dim{\operatorname{dim}}

\def\rk{\operatorname{rk}}
\def\det{\operatorname{det}}

\def\Re{\operatorname{Re}}
\newcommand{\CC}{{\Bbb C}}
\newcommand{\RR}{{\Bbb R}}

\newcommand{\ZZ}{{\Bbb Z}}

\newcommand{\lan}{\langle}
\newcommand{\ran}{\rangle}

\renewcommand{\a}{\alpha}
\renewcommand{\b}{\beta}
\renewcommand{\c}{\gamma}
\renewcommand{\d}{\delta}

\renewcommand{\l}{\lambda}
\newcommand{\m}{\mu}
\renewcommand{\o}{\omega}

\newcommand{\n}{\nu}

\newcommand{\C}{\Gamma}
\newcommand{\D}{\Delta}
\renewcommand{\L}{\Lambda}

\renewcommand{\S}{\Sigma}

\newcommand{\cald}{{\cal D}}

\newcommand{\calf}{{\cal F}}

\newcommand{\calh}{{\cal H}}

\newcommand{\calo}{{\cal O}}

\title{Stationary Phase Integrals, Quantum Toda \\ Lattices, \\
Flag Manifolds and the Mirror Conjecture}

\author{Alexander Givental \thanks{ Research supported by NSF grant 
DMS-93-21915 and
by Alfred P. Sloan Fellowship.} \\ UC Berkeley}

\date{ July 8, 1996,  Revised December 1, 1996}

\begin{document}

\maketitle

{\bf 0. Introduction.} Consider the differential operator
\[ H = \frac{\h ^2}{2} \sum_{i=0}^n \frac{\p ^2}{\p t_i^2} 
- \sum_{i=1}^n e^{t_i-t_{i-1}} \ .\]
This operator is a quantization of the Hamiltonian of the Toda lattice
on $n+1$ identical particles with configuration coordinates $t_0,...,t_n$ 
and with the 
exponential interaction potential $\exp (t_i-t_{i-1})$ of neighbors.
The Toda lattice is known to be integrable on both classical and quantum 
levels:
there exist commuting differential polynomials $D_m(\h \p /\p t,\exp t, \h ),
m=0,...,n$, which play the role of quantum conservation laws (\this 
$[H,D_0]=...=[H,D_n]=0$) and whose symbols $D_m(p,\exp t, 0)$ form a complete
set of Poisson-commuting first integrals of the classical Toda lattice.
In this paper we study solutions $S(t)$ of the differential system 
$D_0 S =...=D_n S =0$ whose characteristic Lagrangian variety $L$ is the most
degenerate invariant Lagrangian variety of the Toda lattice. According to 
\cite{GK, K} this Lagrangian variety is the spectrum of the quantum cohomology 
algebra of the manifold of complete flags $0\subset \CC ^1 \subset ...
\subset \CC ^{n+1}$. We represent solutions $S$ by stationary phase integrals
in $n(n+1)/2$ complex variables and point out the role these solutions 
play in the quantum cohomology theory. As we explain in 
the last section, our results prove the mirror conjecture \cite{Gh} in the
case of the flag manifolds.

\bigskip

{\bf 1. The Toda $\cald$-module.} 

Denote $D_0,...,D_n$ 
coefficients of the polynomial \[ \l^{n+1}+D_0\l^n+...+D_n = \] 
\[ = \det \left[ \begin{array}{ccccc}
\l+p_0 & q_1 &  0  &  0  & ... \\
-1  &\l+p_1 & q_2 &  0  & ... \\
 0  &  -1 &\l+p_3 & q_3 & ... \\
 \  &  .  &  .  &  .  & \   \\
 0  & ... &  0  & -1  &\l+p_n 
\end{array} \right] \ , \]
$D_0=\sum p_i, D_1=\sum_{i> j} p_ip_j + \sum q_i, ...\ $.
The polynomials $D_m(p,q)$ of $p_0,...,p_n$ and \footnote{Throughout 
Sections $1$ -- $3$ we will systematically use the notation $\log q_i$ for 
$t_i-t_{i-1}$.}
$q_1=\exp (t_1-t_0),..., q_n=\exp (t_n-t_{n-1})$ 
form a complete set of Poisson-commuting first integrals of the
Toda lattice (see for instance \cite{GK} ). Their quantizations
\[ D_m(\h \p /\p t_0,...,\h \p /\p t_n , \exp (t_1-t_0), ..., 
\exp (t_n-t_{n-1})), \ m=0,...,n, \]
are defined unambiguously since any monomial in $D_m$ containing $q_i$
contains neither $p_i$ nor $p_{i-1}$.

\medskip

{\bf Theorem 1.} $[H,D_0]=...=[H,D_n]=0$.

{\em Proof.} The commutator 
$[H,\D ]$ of the Hamiltonian
operator $H$ with the above determinant $\D = \l^{n+1} + D_0 \l^n + ... $
vanishes modulo $\h ^2 $ since symbols of $D_m$ Poisson-commute with the 
symbol $\sum p_i^2/2 - \sum q_i $ of $H$. Also $[H,\D ]$ 
does not contain any terms
of order higher than $\h ^2$ since for any $3$ distinct $i,j,k$
we have $\p ^3 q_l /\p t_i \p t_j \p t_k =0$. Computing symbols of the terms 
proportional to $\h ^2$ we find that the contribution of $\D q_i $ is equal to
$-q_i\p^2 \D (p,q)/\p p_i \p p_{i-1}$, and the contribution of 
$\h ^2\p^2 \D /2 \p t_i^2 $ is $q_i \p^2 \D /2\p p_i \p p_{i-1} +
q_{i+1} \p ^2 \D /2\p p_i \p p_{i+1}$. After summation over $i$ all
these contributions cancel out. $\square $

\medskip

We will study solutions of the PDE system 
\[ D_0 S = ... = D_n S = 0 \ (\text{and hence} \ HS=0), \]
\this the solution sheaf of the left module $\cald / \cald (D_0,...,D_n)$
over the algebra $\cald $ of differential operators 
$D(\h \p /\p t, \exp t, \h )$ with Fourier-polynomial coefficients.
We introduce the {\em symbol} $D(p,\exp t, 0)$ of the differential 
polynomial $D$ and call the Lagrangian variety $L\subset T^*(\CC -0)^{n+1}$ 
given by the equations 
$D_0(p,q)=...=D_n(p,q)=0$ the {\em characteristic} Lagrangian variety of 
this $\cald $-module. According to \cite{Ko} the Lagrangian variety $L$ is 
nonsingular. The operator $D_m$ is weighted-homogeneous of degree 
$m+1$ with respect to the grading $\deg \h =1,\ \deg q_i =2, \ \deg t_i =0$
and therefore $L$ is also weighted-homogeneous with weights $\deg p_i =1,\ 
\deg q_i =2$. For generic $t$ the fibers $L\cap T^*_{\exp t}$ of the projection
$L \to (\CC -0)^{n+1}$ consist of $(n+1)!$ distinct simple points. This 
follows from Sard's lemma, and also can 
be deduced by induction from the continued fraction formula for the 
determinant $\D $ which coincides with the numerator of the 
following rational function of $\l $:
\[ \l +p_0 + \frac{q_1  }{\displaystyle \l +p_1 +{\strut 
             \frac{\cdots \hfill }{\displaystyle \cdots +{\strut 
             \frac{q_n }{\displaystyle \l +p_n}}}}} \ . \]
 
\bigskip

{\bf 2. The stationary phase integrals.} 
Consider the following ``$2$ - dimensional
Toda lattice'' with $(n+1)(n+2)/2$ vertices and $n(n+1)$ edges:
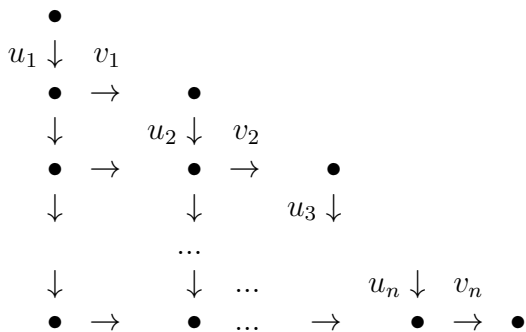
\[ \begin{array}{rrrrrrrr}
 \bullet &  &  &  &  &  &  &  \\
 u_1 \downarrow  & v_1 &  &  &  &  &  &  \\
 \bullet & \rightarrow & \bullet &  &  &  &  &  \\
 \downarrow &   & u_2 \downarrow & v_2 &  &  &  &  \\
 \bullet & \rightarrow & \bullet & \rightarrow & \bullet &   &   &   \\
 \downarrow &  & \downarrow &   & u_3 \downarrow  &   &   &   \\
    &  & ... &  &   &   &   &   \\
 \downarrow &      & \downarrow & ... &     & u_n \downarrow & v_n   &  \\
 \bullet & \rightarrow & \bullet & ... & \rightarrow & 
                                                \bullet & \rightarrow & \bullet
 \end{array} \]

For each edge $\a $ of the lattice we introduce a factor $Q_{\a}$ and
introduce another, excessive notation $u_i$ (respectively $v_i$) for the
factors $Q_{\a}$ corresponding to the vertical (respectively horizontal)
edges next to the diagonal boundary of the lattice as shown on the diagram.
We denote $Y$ the affine algebraic variety in the $n(n+1)$-dimensional
complex space with coordinates $Q_{\a}$ given by the $n(n-1)/2$ equations
$Q_{\a}Q_{\b }=Q_{\c}Q_{\d}$ making the diagram ``commutative'':  
one equation for each $1\times 1$ cell 
\[ \begin{array}{ccccc}
    \   &   \     &    Q_{\a}   &  \      & \  \\
    \   & \bullet & \rightarrow & \bullet & \  \\
 Q_{\c} & \downarrow & \   & \downarrow   & Q_{\b } \\
    \   & \bullet & \rightarrow & \bullet & \  \\
    \   &  \      &     Q_{\d } &   \     & \    \end{array}\ \ . \]
We put $q_1=u_1v_1,...,q_n=u_nv_n$ and denote $Y_q$ the $n(n+1)$-dimensional
fibers of the map from $Y$ to $\CC ^n$ defined by these formulas.

If $q_i\neq 0$ for all $i=1,...,n$, the relations between $Q_{\a}$
allow to express them via the set $\{ T_{\n }, \n=(i,j), 0\leq j\leq
i\leq n \} $ of vertex variables: 
\[ Q_{\a}=\exp (T_{\n _+(\a)} -
T_{\n _-(\a)}),\] 
where $\n _-(\a)$ and $\n_+(\a )$ are respectively the
indices of the source and target vertices of the edge $\a $. In
particular, $q_i=\exp (t_i-t_{i-1})$ where
$t_0=T_{00},t_1=T_{11},...$, and the fiber $Y_q$ is isomorphic to the
$n(n+1)/2$-dimensional complex torus with coordinates $\{ \exp
T_{ij},\ 0\leq j < i \leq n \} $.

On $Y_q$ with $q\in (\CC -0)^n$ we introduce the holomorphic volume form
\[ \o _q = \w _{i=1}^n \w _{j=0}^{i-1} \ d T_{ij} \ ,\]
the holomorphic function $\calf _q$ obtained by restriction to $Y_q$ of
 ``the total Toda potential energy''
\[ \calf = \sum _{\text{edges}\ \a } Q_{\a} \]
and the {\em stationary phase integral}
\[ S_{\C} (t) = \int _{\C \subset Y_q} \ e^{\calf _q /\h } \ \o _q \ .\] 
In this definition, $\C $ represents a (possibly non-compact) cycle in $Y_q$
of middle dimension such that the integral converges. For $\h >0 $ and 
generic $q$ one can construct such a cycle by picking a non-degenerate 
critical point of $\calf _q$ and taking the union of  
descending gradient trajectories of the function $\Re \calf _q$ with respect 
to a suitable Riemannian metric on $Y_q$ on the role of $\C $. 
The stationary phase integral
depends only on the homology class of the cycle in the appropriate homology
group $\calh _q$ which can be described as the inverse limit as 
$M\to \infty $ of the relative homology groups 
\[ H^{\dim _{\CC } Y_q} ( Y_q , \{ \Re \calf _q \leq -M \} ) .\]
The rank of the group $\calh _q$ is equal to the number of critical points
of $ \calf _q$ for generic $q$ since all critical points of the real part of a
holomorphic Morse function have the same Morse index.
   The Gauss-Manin parallel transport of cycles identifies the groups $\calh_q$
for close $q$ but may give rise to a nontrivial global monodromy. The notation
$S_{\C }(t)$ emphasizes the multiple-valued character of the stationary phase 
integrals which therefore depend on the coordinates 
$\log q_i =t_i-t_{i-1}$ on the 
universal covering of the parameter space. The dependence of the integral
on $\h $ is suppressed in this notation.

\medskip

{\bf Theorem 2.} {\em The stationary phase integrals $S_{\C}(t)$ satisfy the
differential equations $D_0 S = D_1 S = ... = D_n S =0$.} 

{\em Proof.} Application of the differential operator 
$\D = \l^{n+1} + D_0 \l^n + ... $ to the stationary phase integral produces
 the amplitude factor $e^{-\calf /\h }\D e^{\calf /\h} = \det (\l + A)$
where
\[ A=\left[ \begin{array}{ccccc}
-u_1  & u_1v_1  & 0       & ...    &      \\
-1    & v_1-u_2 & u_2v_2  & 0      & ...  \\
0     &  -1     & v_2-u_3 & u_3v_3 &  ... \\
      &  .      &   .     &   .    &      \\
      &  ...    &     0   &   -1   &  v_n \\
\end{array} \right] \ .\]     
We should show that this amplitude is congruent to $\l^{n+1}$ modulo linear
combinations of Lie derivatives along vector fields $Q^m  \p /\p T_{\n }$
tangent to $Y_q$, \this linear combinations of $\h \p Q^m /\p T_{ij}+
Q^m \p \calf / \p T_{ij}$ with $j < i $.
Notice that $\p F /\p T_{\n }$ is the sum of $Q_{\a}$ over the ($\leq 2$) 
edges $\a $ ingoing the vertex $\n $ minus such a sum over the outgoing edges.
  
We begin induction on $n$ by noticing that for $n=1$ the amplitude is
equal to $\l ^2 -\l \p \calf /\p T_{21} $ and apply the induction hypothesis
to the triangular lattice with the principle diagonal $i=j$ cut off.

Consider the differential operator $\D '$ defined as the determinant of the
$3$-diagonal $n\times n$-matrix with $\l+\h \p /\p T_{i,i-1}, \ i=1,...,n$, on
the principal diagonal, $-1$'s under the diagonal and $v_i u_{i+1}, \ i=1,...,
n-1$, above the diagonal. Denote $\calf '$ the sum of all $Q_{\a}$ except 
$u_i$'s and $v_i$'s . By the induction hypothesis (and ``commutativity''
$Q_{\a}Q_{\b}=Q_{\c}Q_{\d}$ of the squares next to the diagonal $i=j$) we may
assume that the amplitude factor 
$\exp (-\calf ' /\h )\D ' \exp (\calf' /\h) $ is
congruent to $\l ^n$ modulo Lie derivatives along $Q^m\p /\p T_{ij}$ with
$|i-j|>1$.

The vector field $Q^m \p/\p T_{i,i-1}$, with no edges adjacent to the
vertex $\n =(i,i-1)$ present in the monomial $Q^m$,  produces the amplitude
$Q^m\p \calf /\p T_{i,i-1}$. Since the vertices $(i,i-1)$ 
do not have common 
edges, addition of such amplitudes allows to transform the amplitude factor
$\exp (-\calf '_q/\h )\D '\exp (\calf '_q/\h )$ purely
algebraicly as if $\p \calf /\p T_{i,i-1}=0$. 
Using such transformations we can 
replace $\p \calf '/ \p T_{i,i-1} $ by $v_i-u_i$, and
the induction hypothesis can be reformulated as the congruence 
to $\l ^{n+1}$ of the $(n+1)$-determinant $\det (\l + B)$ where 
\[ B=\left[ \begin{array}{cccccc}
v_1-u_1 & v_1 u_2 & 0       & ...     &         & 0 \\
-1      & v_2-u_2 & v_2u_3  & 0       & ...     & 0 \\
0       & -1      & v_3-u_3 & v_3 u_4 & ...     & 0 \\
        &    .    &   .     &   .     &         & 0 \\
        &  ...    &    0    &  -1     & v_n-u_n & 0 \\ 
        &         & ...     &   0     &  -1     & 0 
 \end{array} \right] \ . \]   
The matrices $A$ and $B$ admit the factorizations $A=UV, \ B=VU$ into the
product of the following triangular matrices:
\[ U=\left[ \begin{array}{cccccc}
u_1 & 0   & ... &     &      &  \\
1   & u_2 & 0   & ... &      &  \\
0   &  1  & u_3 &  0  & ...  &  \\
    &  .  &  .  &  .  &      &  \\
    & ... &  0  &  1  &  u_n & 0 \\
    &     &     &  0  &  1   & 0      \end{array} \right] , \ 
V=\left[ \begin{array}{cccccc}
-1  & v_1 &  0  & ... &      &   \\
0   & -1  & v_2 &  0  & ...  &   \\
0   &  0  & -1  & v_3 &  0   & ... \\
    &  .  &  .  &  .  &      &   \\
    &     & ... &  0  & -1   & v_n \\
    &     & ... &     &  0   & -1     \end{array} \right] \ .\]
Thus $A$ is similar to $B$ and has the same characteristic polynomial.
$\square $

\medskip

The family $\calf _q$ {\em generates} the following Lagrangian variety 
parametrized
by critical points of the functions (and responsible for stationary phase 
approximations to the integrals $S_{\C}(t)$):
\[ \{ (p,q) \ |  \exists y\in Y_q : \ d_y \calf _q =0, 
\ p=\p \calf _q (y)/\p t \} .\]
The identity $\det (\l + A)=\det (\l +B)$ also proves by induction the
following

\medskip

{\bf Corollary 1.} {\em The Lagrangian variety generated by the family 
$\calf _q$ 
coincides with the invariant Lagrangian variety 
$L=\{ (p,q) | D_0(p,q)=...=D_n(p,q)=0\} $ of the Toda lattice.}

\medskip

Notice that the equations $\p \calf /\p T_{\n}=0$ of the critical points and
the notations $p_i=\p \calf /\p t_i$ can be interpreted, in the spirit
of the elementary theory of linear electric circuits, as homological 
boundary conditions
for the $1$-chain $\sum Q_{\a} [\a ]$ on the oriented graph. One can therefore
describe the critical points by the relations 
$Q_{\a}=J_{\phi_+(\a)}-J_{\phi_-(\a)}$ introducing the $2$-chain 
$\sum J_{\phi } [\phi ]$,
a linear combination of the $n(n+1)/2$ clockwise oriented $1\times 1$-cells
of the lattice (the edge $\pm \a $ occurs in the boundary of the
cells $\pm \phi_{\pm }(\a )$, the index $\phi =(i,j)$ 
runs $1\leq j\leq i \leq n$). Additionally, 
$p_0=-J_{11}, p_1=J_{11}-J_{22}, p_2=J_{22}-J_{33}, ... , p_n =J_{nn}$. 
Now cancellations of all $J_{ij}$
with $i>j$ in the total sum $\sum Q_{\a} $ prove

\medskip

{\bf Corollary 2.} {\em The generating function $\calf _q(y_{crit})$ on the
Lagrangian variety $L$ equals $\ \sum _{i=1}^n 2J_{i,i} = 
-np_0+(2-n)p_1+(4-n)p_2+...+ np_n $.} 

{\em Remark.} Corollary $2$ can be also deduced (see \cite{GK}) from 
the weighted homogeneity of $L$. 

\bigskip

{\bf 3. Quantum cohomology of flag manifolds.}  The cohomology algebra
$H^*(F)$ of the flag manifold $F=\{ 0\subset \CC ^1 \subset \CC ^2 \subset
... \subset \CC^{n+1} \}$ is multiplicatively generated by the $1$-st Chern
classes $p_i$ of the tautological line bundles with fibers $\CC ^{i+1}/\CC ^i$.
A complete set of relations between the generators can be written in the form
$(\l + p_0)(\l +p_1)...(\l +p_n)=\l^{n+1}$ equating the 
total Chern class of the sum of the tautological line bundles to that of 
the trivial bundle with the fiber $\CC ^{n+1}$. 

On the Poincare-dual language of intersection indices $\lan \cdot ,\cdot \ran $
the structural constants $\lan ab,c \ran $ of the cohomology algebra count
(with signs) isolated common intersection points of the three cycles $a,b,c$
in general position. 

The quantum cohomology algebra of the flag manifold is defined as a deformation
of the algebra $H^*(F)$ with structural constants $\lan a\circ b,c \ran $ 
counting
isolated holomorphic spheres $(\CC P^1, 0,1,\infty )\to (F, a,b,c)$ passing by
the three marked points through the three cycles. A precise definition can be
based on Gromov's compactness theorem and Kontsevich's concept \cite{Kn} 
of stable holomorphic maps.

Let $\S $ denote a compact connected complex curve with at most double 
singular points. It is rational (\this $H^1(\S , \calo)=0$) if and only if all 
irreducible components of $\S $ are spheres and the
incidence graph of these components is a connected tree. Denote 
$x=(x_1,...,x_k)$ an ordered set of pairwise distinct non-singular marked
points on $\S $. Two holomorphic maps $(\S, x)\to F, \ (\S ', x') \to F$ are
called {\em equivalent} if they are identified by a holomorphic isomorphism
$(\S ,x)\to (\S ',x')$. A holomorphic map $(\S ,x)\to F$ is called {\em stable}
if it does not admit non-trivial infinitesimal automorphisms. For rational $\S$
stability means that each irreducible component mapped to a 
point in $F$ carries at least $3$ marked or singular points. 

The {\em degree} $d$ of the map $\S \to F$ is defined as the total $2$-nd 
homology class of $F$ it represents. According to \cite{Kn, BM}, equivalence 
classes of stable degree-$d$ holomorphic map of rational curves with $k$ 
marked points to the flag manifold $F$ form a compact
complex non-singular orbifold which we denote $F_{k,d}$. If non-empty, it has 
the dimension
$\dim F + \lan -K_F, d \ran +k-3$  
where $-K_F$ is an anti-canonical divisor of $F$. 

Evaluation of maps $(\S ,x)\to F$ at $x_i$ defines the {\em evaluation
maps} $e_i: F_{k,d}\to F,\ i=1,...,k$.

\medskip

{\em Examples.} 1) The moduli space $F_{3,0}$ consists of classes of {\em 
constant} maps $(\CC P^1, 0,1,\infty) \to F$ and thus is isomorphic to $F$. 
The moduli spaces $F_{k,0}$ with $k<3$ are not defined since constant rational
maps with less than $3$ marked points are unstable.

2) Forgetting the space $\CC ^i$ in a flag defines the projection
$F\to F_i$ to the partial flag manifold $F_i$ with fibers isomorphic to 
$\CC P^1$. Denote $a_i\in H_2(F)$ the homology class of the fiber.
In fact any compact holomorphic curve in $F$ of the degree $a_i$ is one
of the fibers. This identifies the moduli space $F_{0,a_i}$ with the
base $F_i$. 

3) The evaluation map $e_1: F_{1,a_i}\to F$ is an isomorphism.

4) The Borel-Weil representation theory for $SL_{n+1}(\CC )$ identifies the 
root lattice of type $A_n$ with the Picard lattice $H^2(F)$ of the flag
manifold. In $H_2(F)$ the classes $a_1,...,a_n$ form the basis of
simple coroots dual to the basis $J_1,...,J_n\in H^2(F)$ of 
{\em fundamental weights} $J_i=-p_0+...-p_i$, \this Chern classes of the 
tautological bundles $\w ^i(\CC ^i)^*$. The simplicial cone spanned 
in $H^2(F, \RR )$ by the fundamental weights is  
the K\"ahler cone of the flag manifold. This implies
that the degree $d$ of any compact holomorphic curve in $F$ is a sum
$d_1a_1+...+d_na_n$ with all $d_i\geq 0$. The sum $\sum_{i>j} p_i-p_j =
2(J_1+...+J_n)$ of positive roots represents the anti-canonical class of
the flag manifold. Thus $\dim F_{3,d}=\dim F + 2(d_1+...+d_n)$.

\medskip

Let us introduce the grading in the algebra $\CC[\L]=\CC [q_1,...,q_n]$ 
of the semigroup
$\L=\{ d=\sum d_ia_i\in H_2(F)| d_i\geq 0\} $ by putting
$\deg q^d =4(d_1+...+d_n)$. We extend $\CC [\L ]$-bilinearly the Poincare 
pairing $\lan A, B\ran =\int _F A\w B$ in the De Rham cohomology of $F$ to 
the graded
$\CC [\L ]$-module $H^*(F,\CC [\L ])$.
Structural constants of the 
$\CC [\L ]$-bilinear quantum multiplication $\circ $ 
on $H^*(F,\CC [L])$ are defined by
\[ \forall A,B,C\in H^*(F)\ \lan A\circ B,C\ran =\sum_{d\in \L} 
    q^d \int _{F_{3,d}} e_1^*(A)\w e_2^*(B)\w e_3^*(C) \ .\]
One can show (see \fin \cite{BM, Ge}) that 
on $H^*(F,\CC [\L ])$ the quantum multiplication
defines the structure of a commutative associative 
graded Frobenius algebra with unity $1\in H^*(F)$ (\this $\lan A\circ 1, C\ran=
\lan A, C\ran $ ).  Modulo $(q_1,...,q_n)$ 
this structure coincides with the cup-product on $H^*(F)$. Recent developments
in symplectic topology and theory of stable maps allow to extend the 
construction of the quantum cohomology algebra to arbitrary compact symplectic
manifolds. 

\medskip

The spectrum of a quantum cohomology algebra can be naturally identified with
certain Lagrangian variety --- the characteristic Lagrangian variety of
the {\em quantum cohomology $\cald $-module}.

For $k>0$ denote $c$ the $1$-st Chern class of the tautological
line bundle over the moduli space $F_{k,d}$ with the fiber at the point
$[(\S ,x)\to F]$ equal to the cotangent line $T^*_{x_1}\S $ to the curve at
the first marked point. Denote $pt=\sum_{i=0}^n p_it_i=\sum (J_i-J_{i+1})t_i$ 
the general $2$-nd cohomology class of $F$ and 
$dt=\sum_{i=1}^n d_i (t_i-t_{i-1})$ the value of
this cohomology class on the homology class $d=\sum d_ia_i$.
For each $A\in H^*(F)$ the vector-function $s_A(t)$ of
$t$ with values in $H^*(F,\CC )$ is defined by
\[ \forall B\in H^*(F)\  \lan s_A(t), B\ran = \lan e^{pt/\h }A, B\ran +
\sum_{d\in \L-0} e^{dt} \int _{F_{2,d}} \frac{e_1^*(e^{pt/\h}A)}{\h - c} 
\w e_2^*(B)\ .\]
By the definition, $s_A$ is a formal power series of 
$q_i=\exp(t_i-t_{i-1})$ with vector-coefficients which are polynomial in 
$\log q_i$ and $\h ^{-1}$.

The vector series $s_A$ satisfy the following linear differential equations 
with periodical coefficients (see \fin \cite{DVV, Db, Ge}):
\[ \h \frac{\p }{\p t_i} s_A = p_i \circ s_A, \ i=0,...,n, \]
where $p_i\circ $ are operators of quantum multiplication by $p_i$. 
In particular, the equations are compatible for any value of $\h $ (\this
$p_i\circ p_j=p_j\circ p_i$ and 
$\p (p_i\circ )/\p t_j =\p (p_j\circ )/\p t_i$),
and the linear space of all solutions to this system coincides with
the space of all vector-functions $s_A$ (of dimension $\rk H^*(F)$ over, say,
$\CC ((\h ^{-1}))\ $).

By definition, the quantum cohomology $\cald $-module is generated by
the scalar functions $S_A(t):=\lan s_A(t), 1 \ran $, \this coincides
with $\cald /I$ where $I=\{ D\in \cald | DS_A=0 \ \forall A\in H^*(F) \} $.
It is easy to show (see for instance \cite{Ge}) that if a homogeneous 
differential polynomial $D(\h \p /\p t, \exp t, \h )$ annihilates 
all the functions $S_A$ then the relation $D(p,q,0)=0$ holds in the 
quantum cohomology algebra.  

\medskip

{\em Examples.} 5) By the very definition $\sum _i \p S_A/\p t_i =0$ since
$\sum p_i =0 $ in $H^*(F)$.

6) The relation $p_0^2+...+p_n^2=0$ holds in the cohomology algebra $H^*(F)$
of the flag manifold. We claim that in the quantum cohomology algebra
$p_0^{\circ 2}+...+p_n^{\circ 2}=2q_1+...+2q_n$. 
Indeed, for the degree reasons 
$J_i\circ J_j$ must be equal to $J_iJ_j$ plus a linear combination 
$\sum \l_m q_m $. The coefficient $\l _m$ is equal to the number
of holomorphic maps $\CC P^1\to F$ passing by $0,1$ and $\infty $ trough
the generic divisors $J_i, J_j$ and a given generic point in $F$. From Examples
$2$ and $4$ we find therefore that $J_i\circ J_j=J_iJ_j+\d _{ij} q_i$ and thus
$\sum p_i^{\circ 2}/2 =\sum J_i^{\circ 2}-\sum J_i\circ J_{i-1} = \sum q_i$.

\medskip

{\bf Theorem 3} $\ $ (\cite{K}). {\em The Hamiltonian operator
$H=(\sum \h ^2\p ^2 /\p t_i^2)/2 - \sum \exp (t_i-t_{i-1})$ 
of the quantum Toda lattice annihilates the functions $S_A, A\in H^*(F)$.}

{\em Proof.} Application of the operator $H$ to $S_A$ yields
\[ \lan H(p\circ ,q)\ s_A, 1\ran + 
\h \sum_i \lan \frac{\p (p_i\circ )}{\p t_i } s_A, 1\ran  \ .\]
The first term vanishes due to Example $6$. Since $\lan  p_i \circ B , 1\ran =
\lan B, p_i \circ 1\ran =
\lan B, p_i \ran $ is constant for any $B\in H^*(F)$, we have
$\lan ( \p p_i \circ /\p t_i) s_A, 1 \ran =
\lan s_A, (\p p_i \circ  /\p t_i) 1\ran =0$. $\ \square $

\medskip

{\bf Corollary.} {\em $\ S_A=\lan e^{pt/\h } s, A\ran $ where the coefficients
$s^{(d)}\in H^*(F,\CC (\h ))$ of the vector-function 
$s=\sum_{d\in \L} s^{(d)} q^d $ can be found recursively from} 
\[ s^{(0)}=1,\ \h [\h (d, d)+\sum_{i=1}^n d_i J_i] s^{(d)}
= \sum_{i: d_i>0} s^{(d-a_i)} \ .\]
{\em In particular, the formal series $s$ converges everywhere.}

\medskip

{\em Remark.} Analogous computations for the flag manifold $G/B$ of a
semi-simple complex Lie group $G$ of rank $n$ give rise to the Hamiltonian 
$H=\sum (a_i,a_j)J_iJ_j/2 - \sum (a_i,a_i) q_i/2 $ of the 
Toda lattice corresponding to the system of simple coroots $a_1,...,a_n$
Langlands - dual to the root system of $G$.  
The proof of the theorem below is a
specialization to the case $G=SL_{n+1}(\CC )$ of the general results by B.Kim 
describing quantum cohomology $\cald $-modules of flag manifolds $G/B$ 
in terms of quantized Toda lattices.

\bigskip

{\bf Theorem 4.} $\ D_0S_A=D_1S_A=...=D_nS_A=0$ {\em for all $A\in H^*(F)$.}

{\em Proof} (see \cite{K}). Since $HD_m S_A=D_mHS_A=0$, 
the formal power series  $S=\sum S^{(d)} q^d := D_mS_A$ with 
coefficients polynomial in $\log q$ (and $\h^{-1}$) satisfies the hypotheses 
of the following 

\medskip

{\bf Kim's Lemma.} {\em If $S^{(0)}=0$ and $HS=0$ then $S=0$.}

{\em Proof} (see \cite{K}) uses only ellipticity of the operator $H$
and polynomiality of its coefficients.  
Consider a non-zero term $S^{(d)}q^d$ of minimal degree and 
pick in it a non-zero monomial term $const\cdot (\log q)^m q^d$ 
of maximal degree $|m|$. Then this term occurs in $HS$ with the coefficient
$(d , d) \cdot const$ which is also non-zero for $d\neq 0$ 
since the symmetric form $(\ ,\ )$ is positively definite. $\ \square $ 

\medskip

{\bf Corollaries.} 
{\em  $(1)$ The quantum cohomology algebra of the flag manifold
$F$ is canonically isomorphic to 
\[ \CC [p_0,...,p_n,q_1,...,q_n]/(D_0(p,q),...,D_n(p,q)) \ .\]

$(2)$ The total multiplicity of critical points of the functions $\calf _q$ 
in the generating family of the characteristic Lagrangian variety $L$
is equal to $\rk H^*(F)=(n+1)!$

$(3)\  \rk \calh _q=\rk H^*(F)$ if all $q_i\neq 0$.

$(4)$ The monodromy representation 
$\ZZ ^{n+1}=\pi _1(\CC -0)^{n+1} \to Aut (\calh _q)$ of the Gauss-Manin 
connection on $\calh _q$ is unipotent and equivalent to the action of 
$\ZZ ^{n+1}$
on $H^*(F,\CC )$ generated by the multiplication operators  
$A\mapsto \exp (2\pi i  p_k),\ k=0,...,n  $.
 
$(5)$ For each $A\in H^*(F)$ the function $S_A$ has
the stationary phase representation 
$\int _{\C \subset Y_q} \exp (F_q /\h) \o _q$
with suitable $\C =\C (A) \in \calh _q \otimes \CC ((\h ^{-1}))$.

$(6)$ Vice versa, the stationary phase integrals $S_{\C }$ admit the series 
expansions $\lan e^{pt/\h } s, A(\C )\ran $.}        

\medskip
{\em Remark.}
As it is shown in \cite{GK}, the Poincare pairing $\lan \ ,\ \ran $ on the 
quantum cohomology algebra $\CC [p,q]/(D_0(p,q),...,D_n(p,q))$ can be described
by the residue formula
\[ \lan A, B\ran (q) = \frac{1}{(2\pi i)^{n+1}} 
\oint \frac{A(p,q)B(p,q)dp_0\w ... \w dp_n}
{D_0(p,q)...D_n(p,q)} \ .\]
Comparison of the stationary phase approximations for the integrals $S_{\C }$
with the formal asymptotics 
\[ S_{A} (t) \sim \h^{\dim F /2} \sum_{p\in L\cap T^*_{\exp t}} C_p(A)
\frac{e^{\sum 2J_i(p)/\h }}{\det ^{1/2} (\p D_i/\p p_j) |_{(p,\exp t)}} \]
for solutions of the quantum cohomology differential equations shows that
the Jacobian $\det (\p D_i/\p p_j)$ computed at generic points $(p,q)\in L$
coincides, up to a constant factor, with the Hessian 
$\det (\p ^2 \calf _q(T)/ \p T_{\n }\p T_{\m })$
of the function $\calf_q $ computed at the non-degenerate critical points 
corresponding to $(p,q)$. It would be interesting to find a direct proof of 
this identity between the two determinants of sizes $n+1$ and $n(n+1)/2$
respectively.  

\bigskip

{\bf 4. The mirror conjecture.} By the mirror conjecture one usually means the
profound equivalence (see \fin \cite{Y, Kh, M}), discovered several years ago
on the basis of string theory, 
between complex and symplectic geometry in Calabi -- Yau
manifolds (\this compact K\"ahler manifolds which admit non-vanishing 
holomorphic volume forms). In particular, the conjecture predicts that the 
quantum cohomology $\cald $-module corresponding to a Calabi -- Yau manifold 
$X$ describes variations of periods of the holomorphic volume form on another 
Calabi -- Yau manifold $Y$, which has the same dimension as $X$ and whose 
Hodge diamond is mirror-symmetric to that of $X$.

In $1993$ we suggested a generalization of the above correspondence
beyond the class of Calabi -- Yau manifolds. \footnote{We should mention 
however that (earlier) E. Witten's paper \cite{W} can be also understood
as a suggestion of a similar generalization. I am thankful to M. Atiyah who
communicated to me this point of view usually ignored by successors of 
\cite{W}.}
Namely, we conjectured (see \cite{Gh}) that
the quantum cohomology $\cald $-module of a compact symplectic $2N$-dimensional
manifold $X$ 
is equivalent to the $\cald $-module generated by stationary
phase integrals $\int _{\C \in Y_q} \exp (\calf _q/\h ) \o _q $, where 
$(Y_q, \calf _q, \o _q)$
is a suitable algebraic family of (possibly non-compact) 
$N$-dimensional complex manifolds $Y_q$, holomorphic functions $\calf _q: 
Y_q \to \CC $ and non-vanishing holomorphic $N$-forms on $Y_q$.

\medskip

The following arguments spoke in favor of such a generalization.

\medskip

$1$) In terms of the characteristic Lagrangian variety 
$L\subset T^*(H_2(X)^{\v} )$ of the quantum cohomology $\cald $-module the
Poincare pairing $\lan \ ,\ \ran $ on the quantum cohomology algebra $\CC [L]$
of the compact symplectic manifold $X$ is given by the formula
\[ \lan A, B\ran (q) = \sum_{p\in T^*_q\cap L} \frac{A(p)B(p)}{\D (p)} \]
(where $\D $ is the restriction to the diagonal in $L\times L$ of the
function representing Poincare-dual class of the diagonal in $X\times X$). 
This formula resembles the residue pairing 
\[ \frac{1}{(2\pi i)^N} \oint a(y)b(y)\frac{dy_1\w ...\w dy_N}
{\frac{\p \calf _q}{\p y_1}...\frac{\p \calf _q}{\p y_N}} =
\sum_{y_*: d_y\calf _q(y_*)=0} \frac{a(y_*)b(y_*)}
{\det (\p ^2\calf _q/\p y_i\p y_j)|_{y_*}} \]
in singularity theory, and in some examples, including complex projective
spaces and Grassmannians, can be indeed replaced by some residue formula.

\medskip

$2$) Solutions $s_A$ to the differential equations $\h q_i \p s_A/\p q_i =
p_i\circ s_A$ arising from the multiplication $\circ $ in the quantum
cohomology algebra of $X$ admit asymptotical approximations
\[ \lan s_A, B\ran \sim \h^{N/2} \sum_{p\in T^*_q\cap L} 
C_p(A) \ e^{-K_X(p)/\h} [\frac{B(p)}{\sqrt \D (p)} + \calo (\h )] \]
resembling the stationary phase approximations 
\[ \int _{\C \subset Y_q} e^{\calf _q(y)/\h } b(y) d^N y \sim  
\h^{N/2} \frac{e^{\calf _q(y_*)/\h } b(y_*) }
{\det ^{1/2} (\p ^2 \calf _q/\p y_i\p y_j) |_{y_*}} \ .\] 

\medskip

$3$) A by-product of our proof \cite{Gt} for toric symplectic manifolds
of Arnold's symplectic fixed point conjecture \cite{A} 
was the following
multiplicative structure in the Floer homology (it is at least morally
equivalent to the quantum multiplication). A compact toric symplectic manifold
$X$ with Picard number $n$ can be obtained as a symplectic reduction 
$X=\CC ^N//T^n$ of the linear space by the subtorus $T^n\subset T^N$ of the
maximal torus on a generic level of the momentum map $\CC ^N\to 
\text{Lie}^* T^N \to \text{Lie}^* T^n$. Denote $(m_{ij})_{i=1}^n \ _{j=1}^N$
the matrix of the projection $\text{Lie }^* T^N\to \text{Lie }^* T^n$. 
In the quantum cohomology algebra of $X$ the classes $Q_1,...,,Q_N$  of 
coordinate hyperplane divisors satisfy the 
multiplicative relations $Q_1^{m_{i1}}...Q_N^{m_{iN}}=q_i, i=1,...,n$, 
and on the other hand --- can be expressed via some basis 
in $H^2(X)$ as $Q_j=p_1 m_{1j}+...+p_n m_{nj}, j=1,...,N$. It is 
easy to see that the latter set of additive relations specifies the critical
set of the function $\calf = Q_1+...+Q_N$ restricted to the $(N-n)$-dimensional
complex torus $Y_q$ defined by the multiplicative relations. We arrive 
to the mirror family $(Y_q, \calf _q, \o _q)$ of the toric manifold $X$
by introducing the holomorphic volume form 
$\o _q=d\log Q_1 \w ...\w d\log Q_N \ /\ d\log q_1\w ...\w d\log q_n $
on the torus $Y_q$.

\medskip

$4$) In the case of $X=\CC P^{N-1}=\CC ^N//T^1$ corresponding
 stationary phase integrals 
\[ \int _{\C \subset \{ Q_1...Q_N=q \} } e^{(Q_1+...+Q_N)/\h }\  
\frac{dQ_1\w ... \w dQ_N}{dq} \]
satisfy {\em the same} differential equation
$(\h q d /dq)^N S = q S$ as the series
\[ s=e^{(p\log q)/\h } \sum_{d=0}^{\infty } 
\frac{q^d}{(p +\h)^N (p+2\h )^N ... (p+d\h )^N} \ \ \text{mod}\ p^N \]
generating the quantum cohomology $\cald $-module of the complex projective
space (see \cite{G1}). This confirms the mirror conjecture for $\CC P^{N-1}$. 

\medskip

$5$) The above mirror family $(Y_q, \calf _q, \o _q)$ of a toric 
symplectic manifold $X$ agrees with Batyrev's mirrors of Calabi -- Yau
anti-canonical hypersurfaces $X'\subset X$. According to \cite{B} one
can construct Calabi -- Yau mirror manifolds $Y'$ as anti-canonical divisors
in the toric manifold obtained by dualization of the momentum polyhedron
for $X$. In fact periods of holomorphic volume forms on Batyrev's
mirrors can be expressed in terms of our mirrors $(Y_q, \calf _q, \o _q)$ as
the formal Laplace transform
\[ \int _{\c \subset Y_q\cap \calf _q^{-1} (1)} \frac{\o _q}{d \calf _q} \]
of the stationary phase integrals. 
\footnote{By the way, this observation suggests
a formulation of a ``quantum Lefschetz 
hyperplane section theorem'' relating the quantum
cohomology $\cald $-module of a Fano manifold $X$ with that of the 
anti-canonical hypersurfaces $X'\subset X$ 
via such a Laplace transform, \this without any mentioning mirror 
manifolds.}

\medskip

In \cite{Ge} we proved the mirror conjecture for Fano ($K_X < 0$) and 
Calabi -- Yau ($K_X=0$) projective complete intersections. Namely, the 
solutions $S_A$ of the quantum cohomology differential equation for $X$ 
given in $\CC P^{N-1}$ by $r$ 
equations of degrees $l_1+...+l_r\leq N$ are described by the following 
integrals $\int _{\C \subset Y_q} \exp (\calf _q/\h)\ \o _q$:
partition the variables $Q_1,...,Q_N$ into $r+1$ groups of lengths 
$l_0+l_1+...+l_r=N$ and denote $G_0,...,G_r$ the total sums $\sum Q_{\a}$ in
each group (for example, $G_0=Q_1+...+Q_{l_0}$, etc.). Then 
\[ Y_q=\{ (Q_1,...,Q_N) \ |\ Q_1...Q_N=q,\ G_1(Q)=...=G_r(Q)=1 \}, \]
\[  \calf _q =G_0 |_{Y_q},\ \ 
 \o _q =\frac{d\log Q_1\w ...\w d\log Q_N}
{d\log q\w d G_1\w ... \w d G_r}\ .\]
These formulas include the traditional Calabi -- Yau mirror phenomenon 
as a degenerate particular case with $l_0=0$ and $G_0=0$.

\medskip

The very idea that the mirror phenomenon exists beyond the class of 
Calabi -- Yau manifolds has not attracted much attention of specialists.
One of the causes is that the results mentioned in this section do
not go further than toric complete intersections. 
It was the actual purpose of the 
project \cite{GK} on flag manifolds started in $1993$ by the author and 
B. Kim to improve this situation and enlarge the supply 
of examples confirming the conjecture. Although the discovered 
relation with Toda lattices has been analyzed in the literature
on quantum cohomology of flag manifolds 
(see \fin \cite{Ko, GFP, K}), their mirrors have escaped us so far.   

In the recent paper \cite{EHX} T. Eguchi, K. Hori and C.-S. Xiong,
independently on our lecture \cite{Gh}, have arrived
 to a similar idea of extending the mirror
conjecture to Fano manifolds. They illustrate the idea with the examples
of complex projective spaces, several rational surfaces, Grassmannians 
$G_{4,2}$ and $G_{5,2}$ and extrapolate formulas from the latter 
examples to general Grassmannians. These formulas served us 
as the starting point
for the present paper; suitably modified, they give rise to a 
construction of mirrors for the flag manifolds. Indeed, the construction of 
$(Y_q, \calf _q, \o _q)$ in Section $2$ 
provides a stationary phase integral representation for the solutions $S_A$ 
of the differential equations defined by the quantum cohomology algebra of
the flag manifold.

I would like to thank M. Kontsevich who brought the paper \cite{EHX} to my
attention and B. Kim for communicating his results \cite{K} based on quantum 
Toda lattices.  I am also thankful to B. Kostant who explained to me that
generalizations of the present paper to arbitrary semi-simple Lie algebras
should intertwine mirror manifolds with Whittaker modules \cite{Kw}.


\begin{thebibliography}{20}

\bibitem{A} V.I. Arnold. {\em Sur un propri\'et\'e topologique des 
applications canoniques de la mecanique classique.} C. R. Acad. Sci. Paris
{\bf 261} (1965), 3719 -- 3722.

\bibitem{B} V.V. Batyrev. {\em Dual polyhedra and mirror symmetry for
Calabi -- Yau hypersurfaces in toric manifolds.} Preprint, alg-geom/9310003.

\bibitem{BM} K. Behrend, Yu. Manin. {\em Stacks of stable maps and Gromov --
Witten invariants.} Preprint, 1995.

\bibitem{DVV} R. Dijkgraaf, E. Verlinde, H. Verlinde. {\em Notes on topological
string theory and $2D$ quantum gravity.} Proc. of the Trieste Spring School, 
1990, M. Green et al., eds., World Sci., Singapore, 1991.

\bibitem{Db} B. Dubrovin. {\em Geometry of 2D topological field theories.}
In: Springer Lecture Notes in Math., {\bf 1620} (1996), 120 -- 348.

\bibitem{EHX} T. Eguchi, K. Hori, C.-S. Xiong. {\em Gravitational quantum
cohomology.} Preprint, 1996.

\bibitem{GFP} S. Fomin, S. Gelfand, A. Postnikov. {\em Quantum Schubert 
polynomials.} Preprint, 1996.

\bibitem{Ge} A. Givental. {\em Equivariant Gromov -- Witten invariants.}
Intern. Math. Res. Notices, 1996, No. 13, 1 -- 63.

\bibitem{Gh} A. Givental. {\em Homological geometry and mirror symmetry.} In:
Proceedings of the International Congress of Mathematicians, 
Z\"urich, 1994, Birkh\"auser, 1995, v. 1, pp. 472 -- 480.

\bibitem{G1} A. Givental. {\em Homological geometry I: projective 
hypersurfaces.} Selecta Math., New Series, {\bf 1} ($1995$), No. 2, 
325 -- 345. 

\bibitem{Gt} A. Givental. {\em A simplectic fixed point theorem for toric
manifolds.} In: The Floer memorial volume. H. Hofer, C. H. Taubes,
A. Weinstein, E. Zehnder (eds.), Progress in Math. {\bf 133}, 
Birkh\"auser, 1995, pp. 445 -- 481.

\bibitem{GK} A. Givental, B. Kim. {\em Quantum cohomology of flag manifolds
and Toda lattices.} Commun. Math. Phys. {\bf 168} (1995), 609 -- 641.

\bibitem{K} B. Kim. {\em Quantum cohomology of flag manifolds $G/B$ and
Toda lattices.} Preprint, alg-geom/9607001.

\bibitem{Kn} M. Kontsevich. {\em Enumeration of rational curves via toric
actions.} In: The moduli space of curves, R. Dijkgraaf, C. Faber, van der Geer
(eds.), Progress in Math. {\bf 129}, Birkh\"auser, 1995, pp. 335 -- 368.

\bibitem{Kh} M. Kontsevich. {\em Homological algebra of mirror symmetry.} In:
Proceedings of the International Congress of Mathematicians, Z\"urich, 1994,
Birkh\"auser, 1995, v. 1, pp. 120 -- 139.

\bibitem{Kw} B. Kostant. {\em On Whittaker vectors and representation theory.}
Invent. Math. {\bf 48} (1978), 101-184.
 
\bibitem{Ko} B. Kostant. {\em Flag manifold quantum cohomology, the Toda 
lattice, and the representation with the highest weight $\rho $.}
Selecta Math. New Series {\bf 2} (1996), 43 -- 91.

\bibitem{M} D. Morrison. {\em Mirror symmetry and moduli spaces.} In:
Proceedings of the International Congress of Mathematicians, Z\"urich, 1994,
Birkh\"auser, 1995, v. 2, pp. 1304 -- 1314.

\bibitem{W} E. Witten. {\em Phases of $N=2$ theories in two dimensions.}
Preprint, IASSNS-HEP-93/3.

\bibitem{Y} S.-T. Yau (ed.). {\em Essays on mirror manifolds.} International
Press, Hong Kong, 1992.

\end{thebibliography}
\enddocument